\documentclass[12pt]{article}
\usepackage{amsmath,amsfonts,amssymb,%amsthm,
amsbsy,euscript}

\newcommand{\cP}{\mathcal{P}}
\newcommand{\cX}{{\EuScript X}}
\newcommand{\bu}{\boldsymbol{u}}
\newcommand{\bx}{\boldsymbol{x}}
\newcommand{\dd}{\partial}

\DeclareMathOperator{\Jac}{Jac}

\title{On the Kontsevich $\star$-\/product associativity mechanism}
\author{R.~Buring\thanks{Johann Bernoulli Institute for Mathematics \& %and
Computer Science, University of Groningen, P.O.~Box~407, 9700~AK Groningen, The~Netherlands.\quad
${}^{\S}$\:%{E-mail}: \texttt{A.V.Kiselev\symbol{"40}rug.nl}.\quad%
Partially supported %in part 
by JBI~RUG project~103511 (Groningen).%
},\quad A.~V.~Kiselev${}^{*,\S}$}
\date{}%{February 29, 2016}

\begin{document}
\maketitle
%%% UDC.
\begin{abstract}\noindent%
The deformation quantization by Kontsevich %(1997) 
is a way to construct an associative non\-com\-mu\-ta\-ti\-ve star\/-\/product $\star=\times+\hbar\,\{\,,\,\}_{\cP}+\bar{o}(\hbar)$ in the algebra of formal power series in~$\hbar$ on a given finite\/-\/dimensional affine Poisson manifold: here 
$\times$~is the usu\-al mul\-ti\-p\-li\-ca\-ti\-on, 
$\{\,,\,\}_{\cP}\neq0$ is the Pois\-son bra\-cket, 
and $\hbar$~is the deformation parameter. The 
%binary operation
pro\-duct~$\star$ is assembled at all powers~$\hbar^{k\geqslant0}$ via %a
summation over a certain %some 
set of weighted graphs 
with $k{+}2$~vertices; for each $k{>}0$, every such graph connects the two co\/-\/mul\-ti\-p\-les of~$\star$ using $k$~copies of~$\{\,,\,\}_{\cP}$. 
  %the Poisson bracket.
Cat\-ta\-neo and Fel\-der %(1999) 
in\-ter\-pre\-t\-ed these topological por\-t\-ra\-its as %the 
genuine Feynman diagrams in the Ikeda\/--\/Izawa model %(1993) 
for quantum gravity.

By expanding the star\/-\/product %s %~$\star$ 
up to~$\bar{o}(\hbar^3)$, i.e.,  with respect to %the 
graphs with at most five vertices but possibly containing loops, we %now 
illustrate the %nontrivial 
mechanism \textsf{Assoc}${}=\diamondsuit\,$(\textsf{Pois\-son}) that converts the Jacobi identity for the bracket%s
~$\{\,,\,\}_{\cP}$ into the associativity of~$\star$.%\\[2pt]
%\textbf{Keywords:} Deformation quantization, associative algebra, Poisson bracket, graph complex, star\/-\/product.
\end{abstract}
 %%% in Russian:

\noindent%
Denote by~$\times$ the multiplication in the commutative associative unital algebra $C^\infty(N^n\to\mathbb{R})$ of scalar functions on a smooth $n$-\/dimensional real manifold~$N^n$. Suppose first that a %possibly, 
non\-com\-mu\-ta\-ti\-ve deformation $\star=\times+O(\hbar)$ of~$\times$ is still unital ($f\star 1=f=1\star f$) and asso\-ci\-a\-ti\-ve, 
%\textsf{Assoc}\,$(f,g,h)\mathrel{{:}{=}}
$(f\star g)\star h = f\star(g\star h)$ for $f,g,h\in C^\infty(N^n)[[\hbar]]$. By taking $3!=6$ copies of the asso\-ci\-a\-ti\-vi\-ty equation for the star\/-\/product~$\star$, we infer that the skew\/-\/symmetric part of the lea\-d\-ing deformation term, $\{f,g\}_{\star}\mathrel{{:}{=}}\smash{\tfrac{1}{\hbar}}\bigl(f\star g-g\star f\bigr){\bigr|}_{\hbar\mathrel{{:}{=}}0}$, is a Poisson bracket.\footnote{%In these terms, 
The left\/-\/hand side of the Jacobi identity $\sum_{\circlearrowright}
\{\{f,g\}_{\star},h\}_{\star}{}=0$ is %would be 
an obstruction to the associativity of the star\/-\/product: whenever the Jacobi identity is violated, one cannot have that %\textsf{Assoc}\,$(f,g,h)=0$
$(f\star g)\star h=f\star(g\star h)$.}

Now %let us have it 
the other way round: can the multiplication~$\times$ on a Poisson manifold~$N^n$ be deformed using the %a given 
bracket~$\{\,,\,\}_{\cP}$ such %so
that the 
$\Bbbk[[\hbar]]$-\/linear star\/-\/product $\star=\times+\hbar\,\{\,,\,\}_{\cP}+\bar{o}(\hbar)$ stays %is
associative? Kontsevich proved~\cite{KontsevichFormality} that on finite\/-\/dimensional affine\footnote{%
%In Remark~\ref{RemAffine} on p.~\pageref{RemAffine} below we recall why the manifold~$N^n$ must be affine, so that 
On affine manifolds~$N^n$,
the only shape %possible form 
of coordinate changes is~$\widetilde{\bu}=A\cdot\bu+\vec{\boldsymbol{c}}$.
  %, cf.~\cite{KontsevichFormality,dq15}. 
%In fact, 
Yet no loss of generality occurs if the space~$N^n$ is the fibre in an affine bundle~$\pi$ of physical %gauge 
fields $\{\bu=\phi(\bx)\}$ over the space\/-\/time%manifold
~$M^m\ni\bx$; % e.g., as in the Yang\/--\/Mills models of the Einstein gravity set\/-\/up. 
the Jacobians $\dd\widetilde{\bu}/\dd\bu=A(\bx)$ are then constant over~$N^n$.
(The arguments of~$\star$ are local functionals of sections, $\phi\in\Gamma(\pi)\to\Bbbk$; 
%The fibre~$N^n$ itself needs not be Poisson;
%%% + The notion of P.br. is then extended.
%for 
the $\star$-\/product is marked by the variational Poisson brackets $\{\,,\,\}_{\boldsymbol{\cP}}$ on the jet space~$J^\infty(\pi)$.)
The deformation quantization from~\cite{KontsevichFormality} is lifted to the gauge field set\/-\/up in~\cite{dq15}.%
}
Poisson manifolds, this is always possible: from~$\{\,,\,\}_{\cP}$ one obtains the bi\/-\/differential terms~$B_k({\cdot},{\cdot})$ at all powers of~$\hbar^{k\geqslant0}$ in the formal series for~$\star$. 
%%%%%%%%
This associative unital $\star$\/-\/product was constructed %by Kontsevich
in~\cite{KontsevichFormality} using a pictorial language:\ the %bi\/-\/differential
operators $B_k%({\cdot},{\cdot})
=\sum\nolimits_{\{\Gamma\}}w(\Gamma)\times
B_k^{\Gamma}({\cdot},{\cdot})$ are encoded by the weighted oriented graphs~$\Gamma$ with $k+2$ vertices and $2k$~edges but without tadpoles or multiple edges;
in every such~$\Gamma$, there are $k$ internal vertices (each of them is a tail for two edges) and $2$ sinks (no issued edges).
 %, see~\cite{KontsevichFormality} or~\cite{dq15} and references therein.
The Poisson bracket $\{\,,\,\}_{\cP}$ with coefficients $\cP^{ij}(\bu)$ at $\bu\in N^n$ provides the ``building block'' ${\pmb\wedge}
=\smash
{\xleftarrow[\text{Left}]{i}\bullet\xrightarrow[\text{Right}]{j}}$ in which 
$\sum_{i,j=1}^n$ is implicit and the vertex contains~$%\hbar
\cP^{ij}(\bu)$.\rule{0pt}{0.9\baselineskip}
To indicate the ordering of indexes in $\cP^{ij}=-\cP^{ji}$,
the out\/-\/going edges are ordered by Left~$\prec$ Right.
%grasping the skew\/-\/symmetry $\cP^{ij}=-\cP^{ji}$.
The edges carry the derivatives $\dd_i\equiv\dd/\dd u^i$ and $\dd_j\equiv\dd/\dd u^j$, respectively. Every such derivation acts on the content of the vertex at the %respective 
arrowhead via the Leibniz rule (and it does so independently from the other in\/-\/coming arrows, if any).\footnote{%We see that 
For example%E.g.
, $%\hbar\,
\{f,g\}_{\cP}(\bu) = 
f\xleftarrow[\text{Left}]{i}\bullet\xrightarrow[\text{R}\smash{\text{ight}}]{j}g =
(f)\overleftarrow{\dd_i}{\bigr|}_{\bu}\cdot %\hbar
\cP^{ij}(\bu)\cdot \overrightarrow{\dd_j}{\bigr|}_{\bu}(g)$, see~\eqref{FigOh3} above.}
  %e.g., as in formula~\eqref{FigOh3} below.

The weights\footnote{Willwacher and Felder (%IMRN 
2010) conjecture that the weights can be %are 
\emph{irrational} numbers for some %all 
graphs.}
$w(\Gamma)\in\mathbb{R}$ of such graphs~$\Gamma$ are given by the integrals over configuration spa\-ces of $k$~distinct points in the hyperbolic plane~$\mathbb{H}^2$ (e.g., in its upper 
half\/-\/plane %Poin\-ca\-r\'e 
model).%
%%%
\footnote{The wedge factors within the integrand in the formula for~$w(\Gamma)$ are copies of the kernel of the singular linear integral operator $(\mathrm{d}*\mathrm{d})^{-1}$ in the hyperbolic geometry of~$\mathbb{H}^2$, see~\cite{CattaneoFelder2000}. Cattaneo and Felder also showed that the $\star$-\/product of two functions $f,g\in C^\infty(N^n\to\mathbb{C})$ amounts to the Feynman path integral calculation of the correlation function,
$\bigl(f\star g\bigr)(\bu)=\smash{\int_{\cX(\infty)=\bu}{ }}\mathrm{D}\cX\,\mathrm{D}\eta\:
f\bigl(\cX(0)\bigr)\times g\bigl(\cX(1)\bigr)\times\exp\bigl(\smash{\tfrac{\boldsymbol{i}}{\hbar}} S\bigl(\cP,[\cX,\eta]\bigr)\bigr)$,
in the Ikeda\/--\/Izawa topological open string model on a disk~$D\simeq%\cong
\mathbb{H}^2$ with boundary~$\dd D\ni 0,1,\infty$; here $\cX\colon D\to N^n$ and $\eta\colon D\to T^*D\otimes\cX^*(T^* N^n)$. All details and further references are found in~\cite{CattaneoFelder2000,Ikeda1994}; still let us remember that within the Ikeda\/--\/Izawa model, the perturbative expansions in~$\hbar$ run, in particular, over the graphs with tadpoles (which must be regularized by hand) but at the same time, those path integral calculations reproduce only the weighted oriented graphs without ``eyes'' (e.g., as in $\gets\cdot{}$\raisebox{-0.5pt}{$\rightleftarrows$}${}\cdot\to$, see Eq.\:\eqref{FigOh3} above). 
%on pp.~\pageref{FigOh3start}--\pageref{FigOh3}).
  %\marginpar{Edit}
Because, to the best of our knowledge, the eye\/-\/containing graphs~$\Gamma_i$ such that $w(\Gamma_i)\neq0$ cannot all at once be eliminated from the star\/-\/product~$\star$ via gauge transformations of its arguments and of its output, see Re\-mark~1 on %footnote
p.~\pageref{FootGaugeOut} and~\cite{KontsevichFormality}, many graphs in the original construction of~$\star$ were not recovered in~\cite{CattaneoFelder2000}. Hence there is an open problem to extend or modify the Ikeda\/--\/Izawa Poisson $\sigma$-\/model such that in the new set\/-\/up, the correlation functions would expand with respect to all the Kontsevich graphs~$\Gamma_i$ with~$w(\Gamma_i)\neq0$.%
}

The associativity postulate for~$\star$ yields the infinite system of quadratic algebraic equations for the weights $w(\Gamma)$ of graphs.%
\footnote{That system solution is not claimed unique: %No uniqueness is claimed for : 
one is provided by the Kontsevich integrals. %Max. generality: balance the gauge classes of graphs.
%%%
Number\/-\/the\-o\-re\-tic %and combinatorial 
properties of those weights were explored by Kontsevich %(1998) 
in the context of motives and %also 
by Will\-wa\-cher\/--\/Fel\-der %(IMRN 2010)
and Garay\/--\/van Straten %(MMJ 2010) 
in the context of Riemann $\zeta$-\/function and Euler $\Gamma$-\/function, respectively.} %This is where number theory meets physics.}
%%%
Kontsevich %proved
shows~\cite{KontsevichFormality} that the left\/-\/hand side 
$\Jac_{\cP}({\cdot},{\cdot},{\cdot})\mathrel{{:}{=}}\sum_{\circlearrowright}\{\{{\cdot},{\cdot}\}_{\cP},{\cdot}\}_{\cP}$ of the Jacobi identity for~$\{\,,\,\}_{\cP}$ is the \emph{only} obstruction to the balance \textsf{Assoc}\,$(f,g,h)\mathrel{{:}{=}}
(f\star g)\star h-f\star(g\star h)=0$ at all powers~$\hbar^k$ of the deformation parameter at once.\footnote{%
Ensuring the associativity \textsf{Assoc}\,$(f,g,h)=0$,
the tri\/-\/vector %Jacobiator 
$\Jac_{\cP}({\cdot},{\cdot},{\cdot})$ is not necessarily (indeed, far not always!\,) evaluated at the three arguments~$f,g,h$ of the associator for~$\star$.}
The core question that we address in this note is how the mechanism 
\textsf{Assoc}${}=\diamondsuit\,$(\textsf{Pois\-son}) works explicitly, making 
%due to which 
the star\/-\/product $\star=\times+\hbar\,\{\,,\,\}_{\cP}+\bar{o}(\hbar)$ associative by virtue of %the 
Jacobi identity for the Poisson bracket~$\{\,,\,\}_{\cP}$.
%%%
Expanding the Kontsevich $\star$-\/product in~$\hbar$ up to~$\bar{o}(\hbar^3)$ and with respect to all the graphs~$\Gamma_i$ such that~$w(\Gamma_i)\neq0$, we obtain\footnote{%N.B.:
%%% M.K.(ord=2, gauge-out the loop at \hbar^2); 
%%% Refer to \ref{RefGauge} where gauging-out the loop at h^2 is explained.
%%% Penkava-Vanhaecke (ord=3, no "eyes").
Balancing the associativity of a star\/-\/product order\/-\/by\/-\/order up to~$\bar{o}(\hbar^3)$, Penkava and Vanhaecke (1998) derived a set of weights for the $(k+2)$-\/vertex Kontsevich graphs \emph{without} loops.
%%% "Ad."
Yet no loops are destroyed %lost
in either of the copies of~$\star$ when the composition $\star\circ\star$ is taken; the associativity of loopless star\/-\/products %$\star^0$ 
is only a part of the full claim for~$\star$.
So, we integrate over the configuration spaces of $k\leqslant3$ points in~$\mathbb{H}^2$ for \emph{all} the Kontsevich graphs (e.g., with loops).%
%%% VVV Repetition? VVV
%those contribute substantially to~$\star$ at all orders~$\hbar^{\geqslant3}$.%
}
%%%%%%%%%
\vspace{-4mm}\label{FigOh3start}
\begin{multline}
\text{\raisebox{-8.5pt}{
\unitlength=0.7mm
\linethickness{0.4pt}
\begin{picture}(12.67,5.67)
\put(2.00,5.00){\circle*{1.33}}
\put(12.00,5.00){\circle*{1.33}}
\put(7.00,5.00){\makebox(0,0)[cc]{$\star$}}
\put(2.00,1.33){\makebox(0,0)[cc]{$f$}}
\put(12.00,1.33){\makebox(0,0)[cc]{$g$}}
\end{picture}
}}
= %\\[-5mm]{=}
\text{\raisebox{-8.5pt}{
\unitlength=0.7mm
\linethickness{0.4pt}
\begin{picture}(15.00,5.67)
\put(0.00,5.00){\line(1,0){15.00}}
\put(2.00,5.00){\circle*{1.33}}
\put(13.00,5.00){\circle*{1.33}}
\put(2.00,1.33){\makebox(0,0)[cc]{$f$}}
\put(13.00,1.33){\makebox(0,0)[cc]{$g$}}
\end{picture}
}}
{+}\frac{\hbar^1}{1!}
\text{\raisebox{-12pt}{
\unitlength=0.7mm
\linethickness{0.4pt}
\begin{picture}(15.00,16.67)
\put(0.00,5.00){\line(1,0){15.00}}
\put(2.00,5.00){\circle*{1.33}}
\put(13.00,5.00){\circle*{1.33}}
\put(2.00,1.33){\makebox(0,0)[cc]{$f$}}
\put(13.00,1.33){\makebox(0,0)[cc]{$g$}}
\put(7.33,16.00){\circle*{1.33}}
\put(7.33,16.00){\vector(-1,-2){5.00}}
\put(7.33,16.00){\vector(1,-2){5.00}}
\end{picture}
}}
{+}\frac{\hbar^2}{2!}
\text{\raisebox{-12pt}{
\unitlength=0.7mm
\linethickness{0.4pt}
\begin{picture}(15.00,20.67)
\put(0.00,5.00){\line(1,0){15.00}}
\put(2.00,5.00){\circle*{1.33}}
\put(13.00,5.00){\circle*{1.33}}
\put(2.00,1.33){\makebox(0,0)[cc]{$f$}}
\put(13.00,1.33){\makebox(0,0)[cc]{$g$}}
\put(7.67,11.67){\circle*{1.33}}
\put(7.67,20.00){\circle*{1.33}}
\put(7.67,20.00){\vector(-1,-3){4.67}}
\put(7.67,20.00){\vector(1,-3){4.67}}
\put(7.67,11.67){\vector(-3,-4){4.00}}
\put(7.67,11.67){\vector(3,-4){4.33}}
\end{picture}
}}
{+}\frac{\hbar^2}{3}{\Biggl(}
\text{\raisebox{-12pt}{
\unitlength=0.7mm
\linethickness{0.4pt}
\begin{picture}(15.00,17.67)
\put(0.00,5.00){\line(1,0){15.00}}
\put(2.00,5.00){\circle*{1.33}}
\put(13.00,5.00){\circle*{1.33}}
\put(2.00,1.33){\makebox(0,0)[cc]{$f$}}
\put(13.00,1.33){\makebox(0,0)[cc]{$g$}}
\put(7.33,11.33){\circle*{1.33}}
\put(2.00,17.00){\circle*{1.33}}
\put(2.00,17.00){\vector(0,-1){11.33}}
\put(2.00,17.00){\vector(1,-1){5.33}}
\put(7.33,11.33){\vector(1,-1){5.33}}
\put(7.33,11.33){\vector(-1,-1){5.33}}
\end{picture}
}}
{+}
\text{\raisebox{-12pt}{
\unitlength=0.7mm
\linethickness{0.4pt}
\begin{picture}(15.00,18.00)
\put(0.00,5.00){\line(1,0){15.00}}
\put(2.00,5.00){\circle*{1.33}}
\put(13.00,5.00){\circle*{1.33}}
\put(2.00,1.33){\makebox(0,0)[cc]{$f$}}
\put(13.00,1.33){\makebox(0,0)[cc]{$g$}}
\put(7.33,11.33){\circle*{1.33}}
\put(7.33,11.33){\vector(1,-1){5.33}}
\put(7.33,11.33){\vector(-1,-1){5.33}}
\put(13.00,17.33){\circle*{1.33}}
\put(13.00,17.33){\vector(0,-1){11.67}}
\put(13.00,17.33){\vector(-1,-1){5.33}}
\end{picture}
}}
{\Biggr)}+\frac{\hbar^2}{6}
\text{\raisebox{-12pt}{
\unitlength=0.7mm
\linethickness{0.4pt}
\begin{picture}(15.00,20.33)
\put(0.00,5.00){\line(1,0){15.00}}
\put(2.00,5.00){\circle*{1.33}}
\put(13.00,5.00){\circle*{1.33}}
\put(2.00,1.33){\makebox(0,0)[cc]{$f$}}
\put(13.00,1.33){\makebox(0,0)[cc]{$g$}}
\put(2.00,15.00){\circle*{1.33}}
\put(13.00,15.00){\circle*{1.33}}
\put(13.00,15.00){\vector(0,-1){9.33}}
\put(2.00,15.00){\vector(0,-1){9.33}}
\bezier{64}(2.00,15.00)(7.00,9.00)(12.67,15.00)
\bezier{60}(13.00,15.00)(7.00,20.33)(2.67,15.00)
\put(11.67,14.00){\vector(1,1){0.67}}
\put(3.33,16.00){\vector(-1,-1){0.67}}
%%%
\put(3,6.5){{\tiny L}}
\put(12,6.5){\llap{\tiny R}}
\put(13.5,15.5){\tiny ``eye''}
\end{picture}
}}
{+}\\
{}+
\frac{\hbar^3}6{\Biggl(}
\text{\raisebox{-12pt}{
\unitlength=0.70mm
\linethickness{0.4pt}
\begin{picture}(15.00,29.00)
\put(0.00,5.00){\line(1,0){15.00}}
\put(2.00,5.00){\circle*{1.33}}
\put(13.00,5.00){\circle*{1.33}}
\put(2.00,1.33){\makebox(0,0)[cc]{$f$}}
\put(13.00,1.33){\makebox(0,0)[cc]{$g$}}
\put(7.67,11.67){\circle*{1.33}}
\put(7.67,20.00){\circle*{1.33}}
\put(7.67,20.00){\vector(-1,-3){4.67}}
\put(7.67,20.00){\vector(1,-3){4.67}}
\put(7.67,11.67){\vector(-3,-4){4.00}}
\put(7.67,11.67){\vector(3,-4){4.33}}
\put(7.67,28.33){\circle*{1.33}}
\put(7.67,28.33){\circle*{1.33}}
\put(7.67,28.33){\vector(-1,-4){5.67}}
\put(7.67,28.33){\vector(1,-4){5.67}}
\end{picture}
}}
{+}
\text{\raisebox{-12pt}{
\unitlength=0.70mm
\linethickness{0.4pt}
\begin{picture}(15.00,23.67)
\put(0.00,5.00){\line(1,0){15.00}}
\put(2.00,5.00){\circle*{1.33}}
\put(13.00,5.00){\circle*{1.33}}
\put(2.00,1.33){\makebox(0,0)[cc]{$f$}}
\put(13.00,1.33){\makebox(0,0)[cc]{$g$}}
\put(2.00,15.00){\circle*{1.33}}
\put(13.00,15.00){\circle*{1.33}}
\put(13.00,15.00){\vector(0,-1){9.33}}
\put(2.00,15.00){\vector(0,-1){9.33}}
\bezier{64}(2.00,15.00)(7.00,9.00)(12.67,15.00)
\bezier{64}(13.00,15.00)(7.00,20.33)(2.67,15.00)
\put(11.67,14.00){\vector(1,1){0.67}}
\put(3.33,16.00){\vector(-1,-1){0.67}}
\put(7.33,23.00){\circle*{1.33}}
\put(7.33,23.00){\vector(-2,-3){4.33}}
\put(7.33,23.00){\vector(2,-3){4.67}}
%%%
\put(3,6.5){{\tiny L}}
\put(12,6.5){\llap{\tiny R}}
\end{picture}
}}
{+}
\text{\raisebox{-12pt}{
\unitlength=0.70mm
\linethickness{0.4pt}
\begin{picture}(15.00,20.33)
\put(0.00,5.00){\line(1,0){15.00}}
\put(2.00,5.00){\circle*{1.33}}
\put(13.00,5.00){\circle*{1.33}}
\put(2.00,1.33){\makebox(0,0)[cc]{$f$}}
\put(13.00,1.33){\makebox(0,0)[cc]{$g$}}
\put(2.00,15.00){\circle*{1.33}}
\put(13.00,15.00){\circle*{1.33}}
\put(13.00,15.00){\vector(0,-1){9.33}}
\put(2.00,15.00){\vector(0,-1){9.33}}
\bezier{64}(2.00,15.00)(7.00,9.00)(12.67,15.00)
\bezier{60}(13.00,15.00)(7.00,20.33)(2.67,15.00)
\put(11.67,14.00){\vector(1,1){0.67}}
\put(3.33,16.00){\vector(-1,-1){0.67}}
\put(7.00,10.33){\circle*{1.33}}
\put(7.00,10.33){\vector(-1,-1){5.33}}
\put(7.00,10.33){\vector(1,-1){5.33}}
%%%
\put(6,12.5){{\tiny R}}
\put(6,18){{\tiny L}}
\end{picture}
}}
{+}
\text{\raisebox{-12pt}{
\unitlength=0.70mm
\linethickness{0.4pt}
\begin{picture}(15.00,23.33)
\put(0.00,5.00){\line(1,0){15.00}}
\put(2.00,5.00){\circle*{1.33}}
\put(13.00,5.00){\circle*{1.33}}
\put(2.00,1.33){\makebox(0,0)[cc]{$f$}}
\put(13.00,1.33){\makebox(0,0)[cc]{$g$}}
\put(7.33,11.33){\circle*{1.33}}
\put(2.00,17.00){\circle*{1.33}}
\put(2.00,17.00){\vector(0,-1){11.33}}
\put(2.00,17.00){\vector(1,-1){5.33}}
\put(7.33,11.33){\vector(1,-1){5.33}}
\put(7.33,11.33){\vector(-1,-1){5.33}}
\put(7.33,22.67){\circle*{1.33}}
\put(7.33,22.67){\vector(-1,-1){5.33}}
\put(7.33,22.67){\vector(0,-1){10.67}}
\end{picture}
}}
{+}
\text{\raisebox{-12pt}{
\unitlength=0.70mm
\linethickness{0.4pt}
\begin{picture}(15.00,22.67)
\put(0.00,5.00){\line(1,0){15.00}}
\put(2.00,5.00){\circle*{1.33}}
\put(13.00,5.00){\circle*{1.33}}
\put(2.00,1.33){\makebox(0,0)[cc]{$f$}}
\put(13.00,1.33){\makebox(0,0)[cc]{$g$}}
\put(7.33,11.33){\circle*{1.33}}
\put(7.33,11.33){\vector(1,-1){5.33}}
\put(7.33,11.33){\vector(-1,-1){5.33}}
\put(13.00,17.33){\circle*{1.33}}
\put(13.00,17.33){\vector(0,-1){11.67}}
\put(13.00,17.33){\vector(-1,-1){5.33}}
\put(7.33,22.00){\circle*{1.33}}
\put(7.33,22.00){\vector(1,-1){5.33}}
\put(7.33,22.00){\vector(0,-1){9.67}}
\end{picture}
}}
%{+}\\
{+}
\text{\raisebox{-12pt}{
\unitlength=0.70mm
\linethickness{0.4pt}
\begin{picture}(15.00,23.33)
\put(0.00,5.00){\line(1,0){15.00}}
\put(2.00,5.00){\circle*{1.33}}
\put(13.00,5.00){\circle*{1.33}}
\put(2.00,1.33){\makebox(0,0)[cc]{$f$}}
\put(13.00,1.33){\makebox(0,0)[cc]{$g$}}
\put(7.33,11.33){\circle*{1.33}}
\put(2.00,17.00){\circle*{1.33}}
\put(2.00,17.00){\vector(0,-1){11.33}}
\put(2.00,17.00){\vector(1,-1){5.33}}
\put(7.33,11.33){\vector(1,-1){5.33}}
\put(7.33,11.33){\vector(-1,-1){5.33}}
\put(2.00,22.67){\circle*{1.33}}
\bezier{64}(2.00,22.00)(-5.00,13.50)(2.33,4.75)
\put(1,6.33){\vector(1,-1){0.67}}
\put(2.00,22.67){\vector(1,-2){5.00}}
\end{picture}
}}
{+}
\text{\raisebox{-12pt}{
\unitlength=0.70mm
\linethickness{0.4pt}
\begin{picture}(20.00,24.67)
\put(0.00,5.00){\line(1,0){15.00}}
\put(2.00,5.00){\circle*{1.33}}
\put(13.00,5.00){\circle*{1.33}}
\put(2.00,1.33){\makebox(0,0)[cc]{$f$}}
\put(13.00,1.33){\makebox(0,0)[cc]{$g$}}
\put(7.33,11.33){\circle*{1.33}}
\put(7.33,11.33){\vector(1,-1){5.33}}
\put(7.33,11.33){\vector(-1,-1){5.33}}
\put(13.00,17.33){\circle*{1.33}}
\put(13.00,17.33){\vector(0,-1){11.67}}
\put(13.00,17.33){\vector(-1,-1){5.33}}
\put(13.00,24.00){\circle*{1.33}}
\put(13.00,24.00){\vector(-1,-2){6.00}}
\bezier{88}(13.00,24.00)(20.00,17.00)(13.67,6.33)
\put(14.00,7.33){\vector(-1,-3){0.33}}
\end{picture}
}}
{\Biggr)}+{}\\
%\end{multline*}
%\begin{equation}
{}+\frac{\hbar^3}3{\Biggl(}
\text{\raisebox{-12pt}{
\unitlength=0.70mm
\linethickness{0.4pt}
\begin{picture}(15.00,17.66)
\put(0.00,5.00){\line(1,0){15.00}}
\put(2.00,5.00){\circle*{1.33}}
\put(13.00,5.00){\circle*{1.33}}
\put(2.00,1.33){\makebox(0,0)[cc]{$f$}}
\put(13.00,1.33){\makebox(0,0)[cc]{$g$}}
\put(7.33,11.33){\circle*{1.33}}
\put(2.00,17.00){\circle*{1.33}}
\put(2.00,17.00){\vector(0,-1){11.33}}
\put(2.00,17.00){\vector(1,-1){5.33}}
\put(7.33,11.33){\vector(1,-1){5.33}}
\put(7.33,11.33){\vector(-1,-1){5.33}}
\put(7.33,8.00){\circle*{1.33}}
\put(7.33,8.00){\vector(-2,-1){5.33}}
\put(7.33,8.00){\vector(2,-1){5.33}}
\end{picture}
}}
{+}
\text{\raisebox{-12pt}{
\unitlength=0.70mm
\linethickness{0.4pt}
\begin{picture}(15.00,17.99)
\put(0.00,5.00){\line(1,0){15.00}}
\put(2.00,5.00){\circle*{1.33}}
\put(13.00,5.00){\circle*{1.33}}
\put(2.00,1.33){\makebox(0,0)[cc]{$f$}}
\put(13.00,1.33){\makebox(0,0)[cc]{$g$}}
\put(7.33,11.33){\circle*{1.33}}
\put(7.33,11.33){\vector(1,-1){5.33}}
\put(7.33,11.33){\vector(-1,-1){5.33}}
\put(13.00,17.33){\circle*{1.33}}
\put(13.00,17.33){\vector(0,-1){11.67}}
\put(13.00,17.33){\vector(-1,-1){5.33}}
\put(7.33,8.33){\circle*{1.33}}
\put(7.33,8.33){\vector(-2,-1){5.33}}
\put(7.33,8.33){\vector(2,-1){5.33}}
\end{picture}
}}
{\Biggr)}
%{+}\\
{+}
\frac{\hbar^3}6{\Biggl(}
\text{\raisebox{-12pt}{
\unitlength=0.70mm
\linethickness{0.4pt}
\begin{picture}(15.00,23.67)
\put(0.00,5.00){\line(1,0){15.00}}
\put(2.00,5.00){\circle*{1.33}}
\put(13.00,5.00){\circle*{1.33}}
\put(2.00,1.33){\makebox(0,0)[cc]{$f$}}
\put(13.00,1.33){\makebox(0,0)[cc]{$g$}}
\put(2.00,15.00){\circle*{1.33}}
\put(13.00,15.00){\circle*{1.33}}
\put(13.00,15.00){\vector(0,-1){9.33}}
\put(2.00,15.00){\vector(0,-1){9.33}}
\bezier{64}(2.00,15.00)(7.00,9.00)(12.67,15.00)
\bezier{60}(13.00,15.00)(7.00,20.33)(2.67,15.00)
\put(11.67,14.00){\vector(1,1){0.67}}
\put(3.67,15.67){\vector(-1,-1){0.67}}
\put(2.00,23.00){\circle*{1.33}}
\bezier{68}(2.00,23.00)(-5,14.00)(1.33,6.00)
\put(2.00,23.00){\vector(3,-2){10.67}}
\put(1,6.67){\vector(1,-2){0.67}}
%%%
\put(3,6.5){{\tiny L}}
\put(12,6.5){\llap{\tiny R}}
\end{picture}
}}
{+}
\text{\raisebox{-12pt}{
\unitlength=0.70mm
\linethickness{0.4pt}
\begin{picture}(19.67,23.00)
\put(0.00,5.00){\line(1,0){15.00}}
\put(2.00,5.00){\circle*{1.33}}
\put(13.00,5.00){\circle*{1.33}}
\put(2.00,1.33){\makebox(0,0)[cc]{$f$}}
\put(13.00,1.33){\makebox(0,0)[cc]{$g$}}
\put(2.00,15.00){\circle*{1.33}}
\put(13.00,15.00){\circle*{1.33}}
\put(13.00,15.00){\vector(0,-1){9.33}}
\put(2.00,15.00){\vector(0,-1){9.33}}
\bezier{64}(2.00,15.00)(7.00,9.00)(12.67,15.00)
\bezier{60}(13.00,15.00)(7.00,20.33)(2.67,15.00)
\put(11.67,14.00){\vector(1,1){0.67}}
\put(3.33,16.00){\vector(-1,-1){0.67}}
\put(13.00,22.33){\circle*{1.33}}
\put(13.00,22.33){\vector(-3,-2){10.33}}
\bezier{84}(13.00,22.33)(19.67,14.33)(13.67,5.67)
\put(14.33,6.67){\vector(-1,-3){0.33}}
%%%
\put(3,6.5){{\tiny L}}
\put(12,6.5){\llap{\tiny R}}
\end{picture}
}}{+}
\text{\raisebox{-12pt}{
\unitlength=0.70mm
\linethickness{0.4pt}
\begin{picture}(15.00,21.33)
\put(0.00,5.00){\line(1,0){15.00}}
\put(2.00,5.00){\circle*{1.33}}
\put(13.00,5.00){\circle*{1.33}}
\put(2.00,1.33){\makebox(0,0)[cc]{$f$}}
\put(13.00,1.33){\makebox(0,0)[cc]{$g$}}
\put(7.33,11.33){\circle*{1.33}}
\put(7.33,11.33){\vector(1,-1){5.33}}
\put(7.33,11.33){\vector(-1,-1){5.33}}
\put(13.00,17.33){\circle*{1.33}}
\put(13.00,17.33){\vector(0,-1){11.67}}
\put(13.00,17.33){\vector(-1,-1){5.33}}
\put(2.00,22.67){\circle*{1.33}}
\put(2.00,22.67){\vector(0,-1){16.67}}
\put(2.00,22.67){\vector(2,-1){10.33}}
\end{picture}
}}
{+}
\text{\raisebox{-12pt}{
\unitlength=0.70mm
\linethickness{0.4pt}
\begin{picture}(15.00,23.33)
\put(0.00,5.00){\line(1,0){15.00}}
\put(2.00,5.00){\circle*{1.33}}
\put(13.00,5.00){\circle*{1.33}}
\put(2.00,1.33){\makebox(0,0)[cc]{$f$}}
\put(13.00,1.33){\makebox(0,0)[cc]{$g$}}
\put(7.33,11.33){\circle*{1.33}}
\put(2.00,17.00){\circle*{1.33}}
\put(2.00,17.00){\vector(0,-1){11.33}}
\put(2.00,17.00){\vector(1,-1){5.33}}
\put(7.33,11.33){\vector(1,-1){5.33}}
\put(7.33,11.33){\vector(-1,-1){5.33}}
\put(13.00,22.67){\circle*{1.33}}
\put(13.00,22.67){\vector(0,-1){16.67}}
\put(13.00,22.67){\vector(-2,-1){10.33}}
\end{picture}
}}
{\Biggr)}{}+\overline{o}(\hbar^3).\label{FigOh3}
%\end{equation}
\end{multline}

%Expanding the associator $\textsf{Assoc}\,(f,g,h) = f \star (g \star h) - (f \star g) \star h$ for $\star$ given by \eqref{FigOh3} up to $\overline{o}(\hbar^3)$, we obtain $6$ terms at $\hbar^1$, $38$ terms $\sim \hbar^2$, and $218$ terms $\sim \hbar^3$.
%Using the antisymmetry $\cP^{ij} = -\cP^{ji}$ to collect the similar terms, we infer that $\textsf{Assoc}\,(f,g,h) = 0$ starts at $\hbar^2$ with $2/3$ times the Jacobi identity,

In every composition $\star \circ \star$ the sums of graphs act on sums of graphs by linearity; %and 
each incoming edge acts via the Leibniz rule (see above).
The mechanism for \textsf{As\-soc}\,$(f,g,h)$ to vanish is two\/-\/step:
first, the sums in $\star \circ \star$ are reduced using the antisymmetry of the Poisson bi\/-\/vector~$\cP$. The output is then reduced modulo the (con\-se\-qu\-e\-n\-c\-es of) Jacobi identity,\footnote{By default, the $L \prec R$ edge ordering equals the left $\prec$ right direction in which edges %they 
start on these pages.}\\[-2mm]
\begin{equation}
\label{eq:Jacobi}
\Jac_\cP(f,g,h) = 
\text{\raisebox{-12pt}{
\unitlength=0.70mm
\linethickness{0.4pt}
\begin{picture}(26.00,16.33)
\put(0.00,5.00){\line(1,0){26.00}}
\put(2.00,5.00){\circle*{1.33}}
\put(13.00,5.00){\circle*{1.33}}
\put(24.00,5.00){\circle*{1.33}}
\put(2.00,1.33){\makebox(0,0)[cc]{$f$}}
\put(13.00,1.33){\makebox(0,0)[cc]{$g$}}
\put(24.00,1.33){\makebox(0,0)[cc]{$h$}}
\put(7.33,11.33){\circle*{1.33}}
\put(7.33,11.33){\vector(1,-1){5.5}}
\put(7.33,11.33){\vector(-1,-1){5.5}}
\put(13,17){\circle*{1.33}}
\put(13,17){\vector(1,-1){11.2}}
\put(13,17){\vector(-1,-1){5.1}}
\end{picture}
}}
{-}
\text{\raisebox{-12pt}{
\unitlength=0.70mm
\linethickness{0.4pt}
\begin{picture}(26.00,16.33)
\put(0.00,5.00){\line(1,0){26.00}}
\put(2.00,5.00){\circle*{1.33}}
\put(13.00,5.00){\circle*{1.33}}
\put(24.00,5.00){\circle*{1.33}}
\put(2.00,1.33){\makebox(0,0)[cc]{$f$}}
\put(13.00,1.33){\makebox(0,0)[cc]{$g$}}
\put(24.00,1.33){\makebox(0,0)[cc]{$h$}}
\put(13,11.33){\circle*{1.33}}
\put(13,11.33){\vector(2,-1){10.8}}
\put(13,11.33){\vector(-2,-1){10.8}}
\put(18.5,17){\circle*{1.33}}
\put(18.5,17){\vector(-1,-1){5.2}}
\put(18.5,17){\vector(-1,-2){5.6}}
\put(13,15){\tiny $L$}
\put(17,12){\tiny $R$}
\end{picture}
}}
{-}
\text{\raisebox{-12pt}{
\unitlength=0.70mm
\linethickness{0.4pt}
\begin{picture}(26.00,16.33)
\put(0.00,5.00){\line(1,0){26.00}}
\put(2.00,5.00){\circle*{1.33}}
\put(13.00,5.00){\circle*{1.33}}
\put(24.00,5.00){\circle*{1.33}}
\put(2.00,1.33){\makebox(0,0)[cc]{$f$}}
\put(13.00,1.33){\makebox(0,0)[cc]{$g$}}
\put(24.00,1.33){\makebox(0,0)[cc]{$h$}}
\put(18.33,11.33){\circle*{1.33}}
\put(18.33,11.33){\vector(1,-1){5.5}}
\put(18.33,11.33){\vector(-1,-1){5.5}}
\put(13,17){\circle*{1.33}}
\put(13,17){\vector(-1,-1){11.2}}
\put(13,17){\vector(1,-1){5.1}}
\end{picture}
}} = 0.
\end{equation}
%Expanding   up to $\overline{o}(\hbar^3)$, we obtain  
For $\star$ given by~\eqref{FigOh3}, the associator 
%$\textsf{Assoc}\,(f,g,h) = (f \star g) \star h - f \star (g \star h)$
contains $6$ terms at $\hbar$, $38$ terms $\sim \hbar^2$, and $218$ terms~$\sim \hbar^3$.
After the use of %Using the antisymmetry 
$\cP^{ij} = -\cP^{ji}$, % to collect the similar terms, 
we infer that \textsf{As\-soc}\,$(f,g,h)$ % =0
starts at $\hbar^2$ with $2/3$ times~\eqref{eq:Jacobi}. %the Jacobi identity,
%The $39$ terms at $\hbar^3$ are the subject of the sequel.
%We immediately recognize $2/3$ times
Next, there are $39$ terms at $\hbar^3$; we now examine how their sum~$\mathsf{A}$ vanishes by virtue of \eqref{eq:Jacobi} and its differential 
consequences.\footnote{%
Within the variational geometry of Poisson field models (cf.~\cite{dq15}),
a tiny leak of the associativity for~$\star$ may occur, if it does at all, only at orders $\hbar^{\geqslant4}$ because at most \emph{one} arrow falls on $\Jac_{\boldsymbol{\cP}}({\cdot},{\cdot},{\cdot})$ in the balance \textsf{As\-soc}\,$(f,g,h)=\bar{o}(\hbar^3)$. But unlike the always vanishing first variation of a homologically trivial functional $\Jac_{\boldsymbol{\cP}}({\cdot},{\cdot},{\cdot})\cong0$, its higher\/-\/order variations can be nonzero.}
Of them, three which are the easiest 
to recognize are\footnote{We use the Einstein summation convention; a sum over all indices is also implicit in the graph notation.}
\\[-2mm]
\begin{equation}
\label{eq:trivialconsequence}
\tfrac{2}{3}\cP^{ij} \Jac_\cP(\dd_i f, \dd_j g, h) = \frac{2}{3}\cdot\Biggl(
\!\!\!\!
\text{\raisebox{-18pt}{
\unitlength=0.70mm
\linethickness{0.4pt}
\begin{picture}(26.00,16.33)
\put(0.00,5.00){\line(1,0){26.00}}
\put(2.00,5.00){\circle*{1.33}}
\put(13.00,5.00){\circle*{1.33}}
\put(24.00,5.00){\circle*{1.33}}
\put(7.33,11.33){\circle*{1.33}}
\put(7.33,11.33){\vector(1,-1){5.5}}
\put(7.33,11.33){\vector(-1,-1){5.5}}
\put(13,17){\circle*{1.33}}
\put(13,17){\vector(1,-1){11.2}}
\put(13,17){\vector(-1,-1){5.1}}
\put(2,17){\circle*{1.33}}
\put(2,17){\vector(0,-1){11.4}}
\put(13,6){\vector(1,-4){0.0}}
\qbezier(2,17)(12,17)(13.00,5.00)
\end{picture}
}}
{-}
\text{\raisebox{-18pt}{
\unitlength=0.70mm
\linethickness{0.4pt}
\begin{picture}(26.00,16.33)
\put(0.00,5.00){\line(1,0){26.00}}
\put(2.00,5.00){\circle*{1.33}}
\put(13.00,5.00){\circle*{1.33}}
\put(24.00,5.00){\circle*{1.33}}
\put(13,11.33){\circle*{1.33}}
\put(13,11.33){\vector(2,-1){10.8}}
\put(13,11.33){\vector(-2,-1){10.8}}
\put(18.5,17){\circle*{1.33}}
\put(18.5,17){\vector(-1,-1){5.2}}
\put(18.5,17){\vector(-1,-2){5.6}}
\put(2,17){\circle*{1.33}}
\put(2,17){\vector(0,-1){11.4}}
\put(13,6){\vector(1,-4){0.0}}
\qbezier(2,17)(12,17)(13.00,5.00)
\put(13,15){\tiny $L$}
\put(17,12){\tiny $R$}
\end{picture}
}}
{-}
\text{\raisebox{-18pt}{
\unitlength=0.70mm
\linethickness{0.4pt}
\begin{picture}(26.00,16.33)
\put(0.00,5.00){\line(1,0){26.00}}
\put(2.00,5.00){\circle*{1.33}}
\put(13.00,5.00){\circle*{1.33}}
\put(24.00,5.00){\circle*{1.33}}
\put(18.33,11.33){\circle*{1.33}}
\put(18.33,11.33){\vector(1,-1){5.5}}
\put(18.33,11.33){\vector(-1,-1){5.5}}
\put(13,17){\circle*{1.33}}
\put(13,17){\vector(-1,-1){11.2}}
\put(13,17){\vector(1,-1){5.1}}
\put(2,17){\circle*{1.33}}
\put(2,17){\vector(0,-1){11.4}}
\put(13,6){\vector(1,-4){0.0}}
\qbezier(2,17)(12,17)(13.00,5.00)
\end{picture}
}}
\!\!
\Biggr) = 0,
\end{equation}
\vskip -1mm
\noindent%
as well as $\tfrac{2}{3}\cP^{ij} \Jac_\cP(f, \dd_i g, \dd_j h) = 0$ and 
$\tfrac{2}{3}\cP^{ij} \Jac_\cP(\dd_i f, g, \dd_j h) = 0$.
So, there remain $30$ terms %,
which vanish via \eqref{eq:Jacobi} in a way more intricate than \eqref{eq:trivialconsequence}.
% (by Kontsevich's theorem) 
It is clear that%Clearly,
\begin{equation}
\label{eq:diffJac}
S_f \mathrel{{:}{=}}
\cP^{ij} \dd_j \Jac_\cP(\dd_i f, g, h) =
\text{\raisebox{-16pt}{
\unitlength=0.70mm
\linethickness{0.4pt}
\begin{picture}(26.00,16.33)
\put(-6, 13){\scriptsize $i$}
\put(0.00,5.00){\line(1,0){26.00}}
\put(2.00,5.00){\circle*{1.33}}
\put(13.00,5.00){\circle*{1.33}}
\put(24.00,5.00){\circle*{1.33}}
\put(7.33,11.33){\circle*{1.33}}
\put(7.33,11.33){\vector(1,-1){5.5}}
\put(7.33,11.33){\vector(-1,-1){5.5}}
\put(13,17){\circle*{1.33}}
\put(13,17){\vector(1,-1){11.2}}
\put(13,17){\vector(-1,-1){5.1}}
\put(13,10.3){\oval(30,20)}
\put(-5.7,19.0){\circle*{1.33}}
\put(-5.7,19.0){\vector(1,-2){6.8}}
\put(-5.7,19.0){\vector(2,-1){5.1}}
\end{picture}
}}
\ 
{-}
\text{\raisebox{-16pt}{
\unitlength=0.70mm
\linethickness{0.4pt}
\begin{picture}(26.00,16.33)
\put(-6, 13){\scriptsize $i$}
\put(0.00,5.00){\line(1,0){26.00}}
\put(2.00,5.00){\circle*{1.33}}
\put(13.00,5.00){\circle*{1.33}}
\put(24.00,5.00){\circle*{1.33}}
\put(13,11.33){\circle*{1.33}}
\put(13,11.33){\vector(2,-1){10.8}}
\put(13,11.33){\vector(-2,-1){10.8}}
\put(18.5,17){\circle*{1.33}}
\put(18.5,17){\vector(-1,-1){5.2}}
\put(18.5,17){\vector(-1,-2){5.6}}
\put(13,10.3){\oval(30,20)}
\put(-5.7,19.0){\circle*{1.33}}
\put(-5.7,19.0){\vector(1,-2){6.8}}
\put(-5.7,19.0){\vector(2,-1){5.1}}
\put(13,15){\tiny $L$}
\put(17,12){\tiny $R$}
\end{picture}
}}
\ 
{-}
\text{\raisebox{-16pt}{
\unitlength=0.70mm
\linethickness{0.4pt}
\begin{picture}(26.00,16.33)
\put(-6, 13){\scriptsize $i$}
\put(0.00,5.00){\line(1,0){26.00}}
\put(2.00,5.00){\circle*{1.33}}
\put(13.00,5.00){\circle*{1.33}}
\put(24.00,5.00){\circle*{1.33}}
\put(18.33,11.33){\circle*{1.33}}
\put(18.33,11.33){\vector(1,-1){5.5}}
\put(18.33,11.33){\vector(-1,-1){5.5}}
\put(13,17){\circle*{1.33}}
\put(13,17){\vector(-1,-1){11.2}}
\put(13,17){\vector(1,-1){5.1}}
\put(13,10.3){\oval(30,20)}
\put(-5.7,19.0){\circle*{1.33}}
\put(-5.7,19.0){\vector(1,-2){6.8}}
\put(-5.7,19.0){\vector(2,-1){5.1}}
\end{picture}
}} = 0.
\end{equation}
\vskip -2mm
\noindent%
Working out the Leibniz rule in~\eqref{eq:diffJac},
%Expanding %the left hand side of %by ,
we collect the graphs according to the number of derivatives falling on
% each argument in the triple
each of $(f,g,h)$.
The edge $\begin{picture}(20.0,5.0)\put(0.0, 2.5){\vector(1,0){20.0}\put(-13,-2.0){$j$}}\end{picture}$ provides the differential orders\footnote{In fact, the double edge to $f$ contributes with zero at $(3,1,1)$ due to 
the skew\/-\/symmetry $\cP^{ij} = -\cP^{ji}$.} 
$(3,1,1)$, $(2,2,1)$, $(2,1,2)$, and $(2,1,1)$ twice. % when acting on the internal vertices
Likewise, we see $(1,1,1)$ in~\eqref{eq:Jacobi} and $(2,2,1)$ in~\eqref{eq:trivialconsequence}.

\smallskip\noindent
{\bf Lemma.} A tri-differential operator $\sum_{|I|,|J|,|K|\geqslant 0} c^{IJK}\ \partial_I \otimes \partial_J \otimes \partial_K$ vanishes identically iff all its coefficients vanish: $c^{IJK} = 0$ for every triple $(I,J,K)$ of multi-indices; here $\partial_L = \partial_1^{\alpha_1} \circ \cdots \circ \partial_n^{\alpha_n}$ for a multi\/-\/index $L = (\alpha_1, \ldots, \alpha_n)$.
% TODO: mention/make clear non-constant coefficients
%\quad %\noindent
Moreover, the sums $\sum_{|I|=i,|J|=j,|K|=k} c^{IJK}\ \partial_I \otimes \partial_J \otimes \partial_K$ are then zero for all $(i,j,k)$; in a va\-ni\-sh\-ing 
sum~$X$ of graphs, we denote by~$X_{ijk}$ its vanishing restriction\footnote{%
For example, %\smallskip\noindent%
   %E.g. 
relation \eqref{eq:trivialconsequence} is the consequence of \eqref{eq:diffJac} at order $(2,2,1)$; restriction of~\eqref{eq:diffJac} to~$(2,1,1)$ yields
\vskip -6mm
%The distribution of the $30$ terms over the differential orders is: $4$ terms at $(1,2,1)$, $9$ terms at $(2,1,1)$, $9$ terms at $(1,1,2)$, and $8$ terms at $(1,1,1)$.
%
%Expanding \eqref{eq:diffJac} and restricting to differential order $(2,2,1)$ forces the incoming arrow to land on the second ground vertex, yielding the trivial consequence \eqref{eq:trivialconsequence}. % we saw before.
%Restricting to differential order $(2,1,1)$ instead forces the incoming arrow to land on one of the two internal vertices:
%
\[%\begin{align*}
\Biggl(
\!\!\!
\text{\raisebox{-20pt}{
\unitlength=0.70mm
\linethickness{0.4pt}
\begin{picture}(26.00,16.33)
\put(0.00,5.00){\line(1,0){26.00}}
\put(2.00,5.00){\circle*{1.33}}
\put(13.00,5.00){\circle*{1.33}}
\put(24.00,5.00){\circle*{1.33}}
\put(7.33,11.33){\circle*{1.33}}
\put(7.33,11.33){\vector(1,-1){5.5}}
\put(7.33,11.33){\vector(-1,-1){5.5}}
\put(13,17){\circle*{1.33}}
\put(13,17){\vector(1,-1){11.2}}
\put(13,17){\vector(-1,-1){5.1}}
\put(0.0,13){\circle*{1.33}}
\put(0.0,13){\vector(4,-1){6.5}}
\put(0.0,13){\vector(1,-4){1.8}}
\put(2,13){\tiny $R$}
\end{picture}
}}
\!
{+} \!\!\! %&
\text{\raisebox{-20pt}{
\unitlength=0.70mm
\linethickness{0.4pt}
\begin{picture}(26.00,16.33)
\put(0.00,5.00){\line(1,0){26.00}}
\put(2.00,5.00){\circle*{1.33}}
\put(13.00,5.00){\circle*{1.33}}
\put(24.00,5.00){\circle*{1.33}}
\put(7.33,11.33){\circle*{1.33}}
\put(7.33,11.33){\vector(1,-1){5.5}}
\put(7.33,11.33){\vector(-1,-1){5.5}}
\put(13,17){\circle*{1.33}}
\put(13,17){\vector(1,-1){11.2}}
\put(13,17){\vector(-1,-1){5.1}}
\put(6.0,19){\circle*{1.33}}
\put(6.0,19){\vector(4,-1){6.5}}
\put(6.0,19){\vector(-1,-3){4.4}}
\end{picture}
}}
\!\!\!
\Biggr)
{-}
\Biggl(
\!\!\!
\text{\raisebox{-20pt}{
\unitlength=0.70mm
\linethickness{0.4pt}
\begin{picture}(26.00,16.33)
\put(0.00,5.00){\line(1,0){26.00}}
\put(2.00,5.00){\circle*{1.33}}
\put(13.00,5.00){\circle*{1.33}}
\put(24.00,5.00){\circle*{1.33}}
\put(13,11.33){\circle*{1.33}}
\put(13,11.33){\vector(2,-1){10.8}}
\put(13,11.33){\vector(-2,-1){10.8}}
\put(18.5,17){\circle*{1.33}}
\put(18.5,17){\vector(-1,-1){5.2}}
\put(18.5,17){\vector(-1,-2){5.6}}
\put(9,19){\circle*{1.33}}
\put(9,19){\vector(4,-1){9.5}}
\put(9,19){\vector(-1,-2){6.6}}
\put(17,12){\tiny $R$}
\end{picture}
}}
\!\!
{+} \!\!\!
\text{\raisebox{-20pt}{
\unitlength=0.70mm
\linethickness{0.4pt}
\begin{picture}(26.00,16.33)
\put(0.00,5.00){\line(1,0){26.00}}
\put(2.00,5.00){\circle*{1.33}}
\put(13.00,5.00){\circle*{1.33}}
\put(24.00,5.00){\circle*{1.33}}
\put(13,11.33){\circle*{1.33}}
\put(13,11.33){\vector(2,-1){10.8}}
\put(13,11.33){\vector(-2,-1){10.8}}
\put(18.5,17){\circle*{1.33}}
\put(18.5,17){\vector(-1,-1){5.2}}
\put(18.5,17){\vector(-1,-2){5.6}}
\put(9,19){\circle*{1.33}}
\put(9,19){\vector(1,-2){3.5}}
\put(9,19){\vector(-1,-2){6.6}}
\put(17,12){\tiny $R$}
\end{picture}
}}
\!\!\!
\Biggr)
{-}
\Biggl(
\!\!\!
\text{\raisebox{-20pt}{
\unitlength=0.70mm
\linethickness{0.4pt}
\begin{picture}(26.00,16.33)
\put(0.00,5.00){\line(1,0){26.00}}
\put(2.00,5.00){\circle*{1.33}}
\put(13.00,5.00){\circle*{1.33}}
\put(24.00,5.00){\circle*{1.33}}
\put(18.33,11.33){\circle*{1.33}}
\put(18.33,11.33){\vector(1,-1){5.5}}
\put(18.33,11.33){\vector(-1,-1){5.5}}
\put(13,17){\circle*{1.33}}
\put(13,17){\vector(-1,-1){11.2}}
\put(13,17){\vector(1,-1){5.1}}
\put(6.0,19){\circle*{1.33}}
\put(6.0,19){\vector(4,-1){6.5}}
\put(6.0,19){\vector(-1,-3){4.4}}
\end{picture}
}}
\!\!
{+}
\text{\raisebox{-20pt}{
\unitlength=0.70mm
\linethickness{0.4pt}
\begin{picture}(26.00,16.33)
\put(0.00,5.00){\line(1,0){26.00}}
\put(2.00,5.00){\circle*{1.33}}
\put(13.00,5.00){\circle*{1.33}}
\put(24.00,5.00){\circle*{1.33}}
\put(18.33,11.33){\circle*{1.33}}
\put(18.33,11.33){\vector(1,-1){5.5}}
\put(18.33,11.33){\vector(-1,-1){5.5}}
\put(13,17){\circle*{1.33}}
\put(13,17){\vector(-1,-1){11.2}}
\put(13,17){\vector(1,-1){5.1}}
\put(-0.5,16){\circle*{1.33}}
\put(-0.5,16){\vector(4,-1){18}}
\put(-0.5,16){\vector(1,-4){2.5}}
\put(4,16){\tiny $R$}
\end{picture}
}}
\!\!\!
\Biggr)
= 0.
\]%\end{align*}
\vskip -4mm
\noindent
Similarly, we have $S_g \mathrel{{:}{=}} \cP^{ij} \dd_j %i 
\Jac_\cP(f, \dd_i %j 
g, h) 
= 0$ and $S_h \mathrel{{:}{=}} \cP^{ij} \dd_j %i 
\Jac_\cP(f, g, \dd_i %j 
h) = 0$.%
}
to a fixed differential order~$(i,j,k)$.

The Poisson bi-vector components $\cP^{ij}$ %of the 
can also serve as arguments of the Jacobiator:\footnote{The three tadpoles produce %the terms 
$\Jac_\cP(\partial_i \cP^{ij}, g,h)\,\partial_j f = 0$, which plays its r\^ole in $\mathsf{A}_{111}$ (see the claim below).}
\begin{equation*}%\label{eq:diffJac2}
I_f \mathrel{{:}{=}} 
\dd_j\bigl(\! \Jac_\cP(\cP^{ij}, g, h) \bigr)\,\dd_i f =
\text{\raisebox{-20pt}{
\unitlength=0.70mm
\linethickness{0.4pt}
\begin{picture}(26.00,16.33)
\put(13.00,5.00){\line(1,0){13.00}}
\put(1,-3){\line(3,2){12.00}}
\put(2.00,5.00){\circle*{1.33}}
\put(13.00,5.00){\circle*{1.33}}
\put(24.00,5.00){\circle*{1.33}}
\put(7.33,11.33){\circle*{1.33}}
\put(7.33,11.33){\vector(1,-1){5.5}}
\put(7.33,11.33){\vector(-1,-1){5.5}}
\put(13,17){\circle*{1.33}}
\put(13,17){\vector(1,-1){11.2}}
\put(13,17){\vector(-1,-1){5.1}}
\put(13,10.3){\oval(30,20)}
\put(2,5.0){\vector(0,-1){7.1}}
\qbezier(2,5)(-8,-2)(-2,9)
\put(-2.6,8.0){\vector(1,2){0.5}}
\put(2,-2){\circle*{1.33}}
\put(-7,3){\scriptsize $j$}
\end{picture}
}}
\ 
{-}
\text{\raisebox{-20pt}{
\unitlength=0.70mm
\linethickness{0.4pt}
\begin{picture}(26.00,16.33)
\put(13.00,5.00){\line(1,0){13.00}}
\put(1,-3){\line(3,2){12.00}}
\put(2.00,5.00){\circle*{1.33}}
\put(13.00,5.00){\circle*{1.33}}
\put(24.00,5.00){\circle*{1.33}}
\put(13,11.33){\circle*{1.33}}
\put(13,11.33){\vector(2,-1){10.8}}
\put(13,11.33){\vector(-2,-1){10.8}}
\put(18.5,17){\circle*{1.33}}
\put(18.5,17){\vector(-1,-1){5.2}}
\put(18.5,17){\vector(-1,-2){5.6}}
\put(13,10.3){\oval(30,20)}
\put(2,5.0){\vector(0,-1){7.1}}
\qbezier(2,5)(-8,-2)(-2,9)
\put(-2.6,8.0){\vector(1,2){0.5}}
\put(2,-2){\circle*{1.33}}
\put(17,12){\tiny $R$}
\put(-7,3){\scriptsize $j$}
\end{picture}
}}
\ 
{-}
\text{\raisebox{-20pt}{
\unitlength=0.70mm
\linethickness{0.4pt}
\begin{picture}(26.00,16.33)
\put(13.00,5.00){\line(1,0){13.00}}
\put(1,-3){\line(3,2){12.00}}
\put(2.00,5.00){\circle*{1.33}}
\put(13.00,5.00){\circle*{1.33}}
\put(24.00,5.00){\circle*{1.33}}
\put(18.33,11.33){\circle*{1.33}}
\put(18.33,11.33){\vector(1,-1){5.5}}
\put(18.33,11.33){\vector(-1,-1){5.5}}
\put(13,17){\circle*{1.33}}
\put(13,17){\vector(-1,-1){11.2}}
\put(13,17){\vector(1,-1){5.1}}
\put(13,10.3){\oval(30,20)}
\put(2,5.0){\vector(0,-1){7.1}}
\qbezier(2,5)(-8,-2)(-2,9)
\put(-2.6,8.0){\vector(1,2){0.5}}
\put(2,-2){\circle*{1.33}}
\put(-7,3){\scriptsize $j$}
\end{picture}
}} = 0.
\end{equation*}
\vskip -1mm

\noindent
Likewise, $I_g \mathrel{{:}{=}} \partial_i \bigl(\Jac_\cP(f, \cP^{ij}, h)\bigr)\,\partial_j g = 0$ and $I_h \mathrel{{:}{=}} \partial_i \bigl(\Jac_\cP(f,g,\cP^{ij})\bigr)\,\partial_j h = 0$.
It is the expansion of $I_f$,\ $I_g$,\ $I_h$ via the Leibniz rule that produces the graphs with ``eyes''.
It also yields an order $(1,1,1)$ differential operator on $(f,g,h)$ which cannot be obtained from \eqref{eq:diffJac}. %??? (2) not (4) ?

\noindent
{\bf Claim.} The sum $\mathsf{A}$ of %All the 
$39$ terms at~$%\sim 
\hbar^3$ in \textsf{Assoc}\,$(f,g,h)$ 
vanishes by virtue of %the 
restriction of 
%~\eqref{eq:diffJac} and~\eqref{eq:diffJac2} 
$S_f,S_g,S_h$ and $I_f,I_g,I_h$
to the %seven 
order% triple
s~$(i,j,k)$ %differential orders 
that %which 
are present 
in~$\mathsf{A}$.
%$(1,1,1)$,\ $(1,2,2)$,\ $(2,2,1)$,\ $(2,1,2)$,\ $(1,1,2)$,\ $(1,2,1)$,\ and~$(2,1,1)$.
%\phantom{M}%\quad
%
Indeed, we have\footnote{By using the symbol $%\smash
{\stackrel{\text{\tiny [$m$]}}{=}}$ we indicate the number~$m$ of terms that are eliminated at each step.}
$\mathsf{A}_{221} \smash{\stackrel{\text{\tiny [3]}}{=}}$ $\tfrac{2}{3} (S_f)_{221}$,
$\mathsf{A}_{122} \smash{\stackrel{\text{\tiny [3]}}{=}} \tfrac{2}{3} (S_g)_{122}$, and
$\mathsf{A}_{212} \smash{\stackrel{\text{\tiny [3]}}{=}} -\tfrac{2}{3} (S_h)_{212}$, see~\eqref{eq:trivialconsequence}. %\marginpar{$S_g$\,?}
Finally, we deduce that
%%%%%%%%%%%%%%%%%%%%%%%%%%%%%%%%%%%%
%\noindent
%By taking $\mathsf{A}_{ijk}$ for these triples and recognizing restrictions of \eqref{eq:diffJac} and \eqref{eq:diffJac2}, we conclude that
%%%%%%%%%%%%%%%%%%%%%%%%%%%%%%%%%%%%
%\begin{align*}
$\mathsf{A}_{111} %&
\smash{\stackrel{\text{\tiny [8]}}{=}} \tfrac{1}{6} (I_f - I_h)_{111}$, %&
%\mathsf{A}_{122} \stackrel{\text{\tiny [3]}}{=} -\tfrac{2}{3} (S_h)_{122}, \qquad
%\mathsf{A}_{221} \stackrel{\text{\tiny [3]}}{=} \tfrac{2}{3} (S_f)_{221}, \\
%\mathsf{A}_{212} &\stackrel{\text{\tiny [3]}}{=} -\tfrac{2}{3} (S_h)_{212}, \qquad
$\mathsf{A}_{112} %&
%\smash
{\stackrel{\text{\tiny [9]}}{=}} \bigl(\tfrac{1}{6}I_f + \tfrac{1}{6}I_g - \tfrac{1}{3}S_h\bigr)%{}
_{112}$, %\\
$\mathsf{A}_{121} %&
%\smash
{\stackrel{\text{\tiny [4]}}{=}} \tfrac{1}{3} (I_f - I_h)_{121}$, and %&
$\mathsf{A}_{211} %&
%\smash
{\stackrel{\text{\tiny [9]}}{=}} \bigl(\tfrac{1}{3}S_f - \tfrac{1}{6}I_g - \tfrac{1}{6}I_h\bigr)%{}
_{211}$.
%\end{align*}
%
The total number of terms which we %are 
thus eliminate %d 
equals~$(3 + 3 + 3) + 8 + 9 + 4 + 9 = 39$.\hfill {\footnotesize$\square$}
%, where $3 + 3 + 3$ terms come from \eqref{eq:trivialconsequence}.

\smallskip\noindent%
\textbf{Remark}~1\textbf{.}\quad
The deformation quantization is a gauge theory: each argument~$\bullet$ of~$\star$ marks its gauge %own 
class~$[\bullet]$ under the linear maps %mappings %transformations%\footnote{%
%In fact, the differential operators near $I^{??}$ vanish identically due to~$\cP^{ij}=-\cP^{ji}$.}
%%%
~$\mathsf{t}\colon\bullet\mapsto[\bullet]=\bullet
+\hbar\bigl( I^{\emptyset}\,\dd_i\dd_j(\cP^{ij})^{\equiv0}\times\bullet
+ I^{\circlearrowright}\,\dd_i\cP^{ij}\,\dd_j(\bullet) \bigr)
+{}$\\[-3.5mm]
%%%%%%%%%%%%%%%%%%%%%%%%%%%%% gauge   %%%%%%%%%%%%%%%
\[%\begin{multline*}
\hbar^2\!\!\stackrel{{\scriptscriptstyle(0)}}{I} %C_0
\!\!{\text{\raisebox{-22pt}{
\unitlength=0.7mm
\linethickness{0.4pt}
\begin{picture}(11.67,19.67)
\put(6.00,5.00){\circle*{2.33}}
\put(1.00,15.00){\circle*{1.33}}
\put(11.00,15.00){\circle*{1.33}}
\put(11.00,15.00){\vector(-1,-2){4.67}}
\put(1.00,15.00){\vector(1,-2){4.67}}
\bezier{48}(1.00,15.00)(6.00,11.00)(10.00,14.67)
\bezier{52}(11.00,15.00)(6.00,19.67)(1.33,15.33)
\put(9.33,14.00){\vector(3,2){1.00}}
\put(2.00,16.00){\vector(-3,-2){1.00}}
\end{picture}
}}}
+\hbar^3{\Biggl[%(
}\!\!\stackrel{{\scriptscriptstyle(1)}}{I} %C_1
\!\!{\text{\raisebox{-22pt}{
\unitlength=0.7mm
\linethickness{0.4pt}
\begin{picture}(12.33,21.33)
\put(6.00,5.00){\circle*{2.33}}
\put(6.00,10.00){\circle*{1.33}}
\put(11.33,16.00){\circle*{1.33}}
\put(11.67,17.00){{\tiny${\equiv}0$}}
%%%
\put(0.67,16.00){\circle*{1.33}}
\put(11.33,16.00){\vector(-1,0){10.00}}
\put(0.67,16.00){\vector(3,-4){4.33}}
\put(5.67,10.00){\vector(1,1){5.33}}
\bezier{44}(11.33,16.00)(12.33,10.67)(7.00,10.00)
\bezier{44}(6.00,10.00)(0.33,9.67)(0.67,15.00)
\bezier{56}(0.67,16.00)(6.00,21.33)(10.67,16.67)
\put(9.67,17.67){\vector(3,-2){1.00}}
\put(8.67,10.33){\vector(-4,-1){1.33}}
\put(0.67,13.33){\vector(0,1){1.33}}
\end{picture}
}}}
{+}\!\stackrel{{\scriptscriptstyle(2)}}{I} %C_2
\!\!{\text{\raisebox{-22pt}{
\unitlength=0.7mm
\linethickness{0.4pt}
\begin{picture}(12.33,21.33)
\put(6.00,5.00){\circle*{2.33}}
\put(6.00,10.00){\circle*{1.33}}
\put(11.33,16.00){\circle*{1.33}}
\put(0.67,16.00){\circle*{1.33}}
\put(11.33,16.00){\vector(-1,0){10.00}}
\put(0.67,16.00){\vector(3,-4){4.33}}
\put(5.67,10.00){\vector(1,1){5.33}}
\bezier{44}(11.33,16.00)(12.33,10.67)(7.00,10.00)
\bezier{44}(0.67,16.00)(6,22.67)(10.67,16.33)
\put(9.67,17.67){\vector(3,-2){1.00}}
\put(8.67,10.33){\vector(-4,-1){1.33}}
\put(6.00,10.33){\vector(0,-1){4.67}}
\end{picture}
}}}
{+}\!\stackrel{{\scriptscriptstyle(3)}}{I} %C_3
\!\!{\text{\raisebox{-22pt}{
\unitlength=0.7mm
\linethickness{0.4pt}
\begin{picture}(12.00,22.67)
\put(6.00,5.00){\circle*{2.33}}
\put(11.00,15.00){\circle*{1.33}}
\put(11.33,17.00){{\tiny${\equiv}0$}}
%%%
\put(1.00,15.00){\circle*{1.33}}
\put(6.00,22.00){\circle*{1.33}}
\put(11.00,15.00){\vector(-1,-2){4.67}}
\put(1.00,15.00){\vector(1,-2){4.67}}
\put(6.00,22.00){\vector(-2,-3){4.33}}
\put(6.00,22.00){\vector(3,-4){4.67}}
\bezier{44}(10.67,15.33)(11.00,20.67)(6.67,22.00)
\bezier{40}(1.00,15.00)(0.00,19.33)(5.33,22.00)
\put(4.33,21.67){\vector(2,1){0.67}}
\put(7.00,21.67){\vector(-2,1){0.67}}
%
%\put(55,15){Fig.\:3}
\end{picture}
}}}
{+}\!%\\{+}
\stackrel{{\scriptscriptstyle(4)}}{I} %C_4
\!\!{\text{\raisebox{-22pt}{
\unitlength=0.7mm
\linethickness{0.4pt}
\begin{picture}(12.00,22.66)
\put(6.00,5.00){\circle*{2.33}}
\put(11.00,15.00){\circle*{1.33}}
\put(1.00,15.00){\circle*{1.33}}
\put(6.00,22.00){\circle*{1.33}}
\put(11.00,15.00){\vector(-1,-2){4.67}}
\put(1.00,15.00){\vector(1,-2){4.67}}
\put(6.00,22.00){\vector(-2,-3){4.33}}
\put(6.00,22.00){\vector(3,-4){4.67}}
\put(4.33,21.67){\vector(2,1){0.67}}
\put(11.00,15.00){\vector(-1,0){9.00}}
\bezier{44}(1.00,15.00)(0.00,21.00)(5.00,21.67)
\end{picture}
}}}
{+}\!\stackrel{{\scriptscriptstyle(5)}}{I} %C_5
\!\!{\text{\raisebox{-22pt}{
\unitlength=0.7mm
\linethickness{0.4pt}
\begin{picture}(11.66,23.00)
\put(6.00,5.00){\circle*{2.33}}
\put(1.00,15.00){\circle*{1.33}}
\put(11.00,15.00){\circle*{1.33}}
\put(11.33,17.00){{\tiny${\equiv}0$}}
\put(11.00,15.00){\vector(-1,-2){4.67}}
\put(1.00,15.00){\vector(1,-2){4.67}}
\bezier{48}(1.00,15.00)(6.00,11.00)(10.00,14.67)
\bezier{52}(11.00,15.00)(6.00,19.67)(1.33,15.33)
\put(9.33,14.00){\vector(3,2){1.00}}
\put(2.00,16.00){\vector(-3,-2){1.00}}
\put(6.00,22.33){\circle*{1.33}}
\put(6.00,22.33){\vector(-3,-4){4.67}}
\put(6.00,22.33){\vector(3,-4){5.00}}
\end{picture}
}}}
{+}\!\stackrel{{\scriptscriptstyle(6)}}{I} %C_6
\!\!{\text{\raisebox{-22pt}{
\unitlength=0.7mm
\linethickness{0.4pt}
\begin{picture}(12.67,20.67)
\put(6.00,5.00){\circle*{2.33}}
\put(6.00,15.00){\circle*{1.33}}
\put(12.00,15.00){\circle*{1.33}}
\put(0.00,20.00){\circle*{1.33}}
\put(0.00,20.00){\vector(4,-3){5.33}}
\put(6.00,15.00){\vector(0,-1){9.67}}
\put(0.00,20.00){\vector(1,-3){5}}
\put(12.00,15.00){\vector(-1,-2){5}}
\bezier{36}(6.00,15.00)(8.33,18.67)(11.33,15.33)
\bezier{28}(12.00,15.00)(8.00,12.67)(6.67,14.67)
\put(10.67,16.00){\vector(2,-1){0.67}}
\put(7.33,14.00){\vector(-2,1){0.67}}
\end{picture}
}}}
{+}\!\stackrel{{\scriptscriptstyle(7)}}{I} %C_7
\!\!{\text{\raisebox{-22pt}{
\unitlength=0.7mm
\linethickness{0.4pt}
\begin{picture}(12.33,20.67)
\put(6.00,5.00){\circle*{2.33}}
\put(6.00,20.00){\circle*{1.33}}
\put(0.33,12.33){\circle*{1.33}}
\put(11.67,12.33){\circle*{1.33}}
\put(6.00,20.00){\vector(0,-1){14.33}}
\put(0.33,12.33){\vector(2,3){5.67}}
\put(6.00,20.00){\vector(3,-4){5.33}}
\put(11.67,12.33){\vector(-2,-3){4.67}}
\put(0.33,12.33){\vector(2,-3){4.67}}
\put(5.33,12.33){\vector(-1,0){4.33}}
\put(11.67,12.33){\line(-1,0){5.00}}
\end{picture}
}}}
{\Biggr]%)
}
\]%\end{multline*}
\vskip -3mm
\noindent%
$+\overline{o}(\hbar^3)$,
%%%%%%%%%%%%%%%%%%%%%%%%%%%%%  gauge  %%%%%%%%%%%%%%%%%%%%
where the constants~$\smash{\stackrel{\text{\tiny$(\alpha)$}}{I}}\in\Bbbk$ can be arbitrary\footnote{%COATING:
The view~\cite{CattaneoFelder2000} on $\star$-\/products as $\hbar$-\/expansions of path integrals %for correlators in the model from~\cite{Ikeda94}
shows that the graphs~$\Gamma_i$ in~\eqref{FigOh3} are genuine Feynman diagrams for the channel marked by~$\cP$. %the Poisson bi-vector
The weights $w(\Gamma_i)$ integrate over the energy of each intermediate vertex. %~$\hbar\cP$.
Quite naturally, a particle~$\bullet$ shares its energy\/-\/mass with the interaction carriers~$\cP$ as it gets coated by them. But no object~$\bullet$ can spend more energy on growing its gauge tail than the amount it actually has; hence every set~$[\bullet]$ is bounded in the space of parameters~$\boldsymbol{I}$.%
%(The coating of relativistic particles is a key idea in the Feynman diagram technique.)% by R.~Feynman.
} 
and $\mathsf{t}$~is formally invertible %in
over~$\Bbbk[[\hbar]]$. In turn,
the star\/-\/products are gauged\footnote{\label{FootGaugeOut}%
For example, the loop graph at ${\hbar^2}/{6}$ in~\eqref{FigOh3} is gauged out by~$\mathsf{t}(\bullet)=\bullet+\tfrac{\hbar^2}{12}\!\!  
{\text{\raisebox{-8pt}[1pt][1pt]{%
\unitlength=0.35mm%0.7mm
\linethickness{0.4pt}
\begin{picture}(11.67,19.67)
\put(6.00,5.00){\circle*{2.33}}
\put(1.00,15.00){\circle*{1.33}}
\put(11.00,15.00){\circle*{1.33}}
\put(11.00,15.00){\vector(-1,-2){4.67}}
\qbezier[25](11,15)(6.33,5.67)(6.33,5.67)
\put(1.00,15.00){\vector(1,-2){4.67}}
\qbezier[25](1,15)(5.67,5.67)(5.67,5.67)
%%%
\bezier{48}(1.00,15.00)(6.00,11.00)(10.00,14.67)
\bezier{52}(11.00,15.00)(6.00,19.67)(1.33,15.33)
\put(9.33,14.00){\vector(3,2){1.00}}
\put(2.00,16.00){\vector(-3,-2){1.00}}
\end{picture}
}%raisebox
}%text
}$, see~\cite{KontsevichFormality} for further details.%
}
by using~$\mathsf{t}$: $f\mathbin{{\star}'}g\mathrel{{:}{=}}\smash{\mathsf{t}^{-1}\bigl(\mathsf{t}(f)\star\mathsf{t}(g)\bigr)}$. This
degree of freedom
extends the uniqueness problem for Kontsevich's solution~$\star$ of \textsf{Assoc}\,$(f,g,h)=0$. Namely, not the exact balance of power series but an equivalence~\raisebox{0.67pt}{${\scriptstyle[}$}%
${=}$%
\raisebox{0.67pt}{${\scriptstyle]}$} 
of gauge classes (up to unrelated transformations at all steps) can be sought in $\bigl[ [f]\star[g]\bigr]\star[h]$% [{=}] 
\raisebox{0.67pt}{$\,{\scriptstyle[}$}%
${=}$%
\raisebox{0.67pt}{${\scriptstyle]}\,$}%
$[f]\star\bigl[ [g]\star[h] \bigr]$.

\smallskip\noindent%
\textbf{Remark}~2\textbf{.}\quad
Each graph~$\Gamma$ in~\eqref{FigOh3} encodes the polydifferential operator of scalar arguments in a coordinate\/-\/free way.
The Jacobians $\dd\bu/\dd\widetilde{\bu}$ of \emph{affine} mappings %first
appear on the edges %of~$\Gamma$
but then they join the content~$\hbar\cP^{ij}$ of internal vertices at the arrowtails,\footnote{%For example, 
E.g., 
$\smash{\overleftarrow{\dd'_\alpha}}%\cdot
\widetilde{\cP}{}^{\alpha\beta}
{\bigr|}_{\widetilde{\bu}} %\cdot
\smash{\overrightarrow{\dd'_\beta}} =
\smash{\overleftarrow{\dd_i}}\,\frac{\dd u^i}{\dd\widetilde{u}^\alpha}\,\widetilde{\cP}{}^{\alpha\beta}
{\bigl|}_{\widetilde{\bu}(\bu)} \,\frac{\dd u^j}{\dd\widetilde{u}^\beta}\, \smash{\overrightarrow{\dd_j}} =
\smash{\overleftarrow{\dd_i}}\cdot\cP^{ij}{\bigl|}_{\bu} \cdot
\smash{\overrightarrow{\dd_j}}$ %,
so that
$\{f,g\}_{\cP(\bu)}(\bu)=\{f,g\}_{\widetilde{\cP}(\widetilde{\bu}(\bu))}\bigl(\widetilde{\bu}(\bu)\bigr)$.%
}
forming~$%\hbar
\smash{%\widetilde
\tilde{\cP}}{}^{\alpha\beta}$ from~$%\hbar
\cP^{ij}$.
Independent from~$\bu\in N^n$, these Jacobians stay invisible to all in\/-\/coming arrows (if any). So, the operator given by a graph~$\Gamma$ with~$\hbar\cP(\bu)$ in its vertices is equal to the one for 
$\hbar\smash{%\widetilde
\tilde{\cP}}\bigl(\widetilde{\bu}(\bu)\bigr)$ there.
   %%% vs Gel'fand -- Kazhdan
%%%
This reasoning works for the variational Poisson brackets~$\{\,,\,\}_{\boldsymbol{\cP}}$ on~$J^\infty(\pi)$ for affine bundles~$\pi$ with fibre~$N^n$ over points~$\bx\in M^m$, see~\cite{dq15}. The graphs~$\Gamma$ then yield local variational polydifferential operators 
%%%
%In the paper~\cite{dq15} the second author lifted the Kontsevich construction of star\/-\/products~$\star$ to the %variational 
%Poisson field model geometry, so that the contributions $\hbar^k\,B_k^{\Gamma}({\cdot},{\cdot})$ from graphs~$\Gamma$ are local variational polydifferential operators yet the pictorial language stays intact.
   %In the variational set-up there can be a tiny leak of the associativity at orders \hbar^{\geqslant 3}.
%%%
yet the %pictorial 
pictorial language %of pictures 
%from
of~\cite{KontsevichFormality} is %stays 
the same.%intact.
\footnote{%
A sought\/-\/for %proper 
extension of the Ikeda\/--\/Izawa topological open string geometry
--\,namely, its lift from the %affine 
Poisson manifolds $\bigl(N^n,\{\,,\,\}_{\cP}\bigr)$ in~\cite{CattaneoFelder2000,Ikeda1994} to the variational %Poisson 
set\/-\/up $\bigl(J^\infty(\pi),\{\,,\,\}_{\boldsymbol{\cP}}\bigr)$ of jet spaces in~\cite{dq15}\,-- is a mechanism to quantize Poisson field models. %We shall consider it elsewhere.%
This will be the object of another paper.
}

{\small\smallskip%\paragraph*
\noindent\textbf{Acknowledgements.}
%The second author 
A.\,V.\,K. thanks the organizers %Organizing committee
of international workshop SQS'15 %on {Supersymmetry \& %and 
%Quantum Symmetries}
(Au\-g\-ust 3--8, 2015 at %LTPh 
JINR~Dubna, Russia)
%for a welcome and warm atmosphere during the meeting.
for stimulating %helpful 
discussions and partial financial support.%\\[-2mm] %remarks. %\quad 

%This research was supported in part by JBI~RUG project~103511 (Groningen). 
%A~part of this research was done while the second author was visiting at 
%the $\smash{\text{IH\'ES}}$ (Bures\/-\/sur\/-\/Yvette, France); 
%the financial support and hospitality of this institution are gratefully acknowledged.
}

{\footnotesize%
}

\end{document}